\documentclass[12pt]{amsart}
\usepackage{verbatim}
\usepackage{graphicx}
\usepackage{amssymb}
\usepackage{amsfonts}
\usepackage{amsmath}
\usepackage{hyperref}

\def\R{\mathbb{R}}
\def\RR{\mathbb{R}}

\def\NN{\mathbb{N}}

\def\eps{\varepsilon}

\def\tr{\mathrm{Tr}}

\def\A{\mathcal A}
\def\B{\mathcal B}
\def\sym{\mathrm{Sym}}
\def\Ah{\sym(\A)}
\def\Bh{\sym{\B}}
\def\R{\mathcal R}
\def\sch{\mathrm{sch}}
\def\hmp{\mathrm{hmp}}

\def\rall{\R_\mathrm{all}}

\def\rbd{\R_\mathrm{bnd}}
\def\rfin{\R_\mathrm{fin}}

\def\rev{\R_\mathrm{eval}}
\def\pts{\mathrm{pt}(\R)}

\def\ev{\mathrm{ev}}
\def\T{T}
\def\X{X}

\def\sat{\mathrm{Sat}}

\def\q{e}
\def\fd{F_d}
\def\fdk{F_{d,K}}
\def\pd{\mathcal{P}_d}

\newtheorem{theorem}{Theorem}[section]
\newtheorem{lemma}[theorem]{Lemma}
\newtheorem{thm}[theorem]{Theorem}

\newtheorem{prop}[theorem]{Proposition}

\theoremstyle{definition}

\newtheorem{ex}[theorem]{Example}

\begin{document}
\title[Noncommutative Positivstellens\" atze]{Noncommutative Positivstellens\" atze for pairs representation-vector}
\author{Jaka Cimpri\v c}

\date{July 20th 2010, revised October 5th 2010}

\begin{abstract}
We study non-commutative real algebraic geometry for
a unital associative $\ast$-algebra $\A$ viewing
the points as pairs $(\pi,v)$ where $\pi$ is 
an unbounded $\ast$-representation of $\mathcal A$ 
on an inner product  space which contains the vector $v$.
We first consider the $\ast$-algebras of matrices of usual and free real multivariate polynomials
with their natural subsets of points. If all points are allowed then 
we can obtain results for general $\A$.
Finally, we compare our results with their analogues in the usual
(i.e. Schm\" udgen's) non-commutative real algebraic geometry  
where the points are unbounded $\ast$-representation of $\A$.
\end{abstract}



\address{Cimpri\v c Jakob, University of Ljubljana, Faculty of Mathematics and Physics,
Department of Mathematics, Jadranska 19, SI-1000 Ljubljana, Slovenia,
cimpric@fmf.uni-lj.si,http://www.fmf.uni-lj.si/\~{}cimpric}

\maketitle

\section{Introduction}

Classical real algebraic geometry is interested in \textit{positivity sets} of real multivariate polynomials,
i.e. sets of the form 
\[
K_S=\{a \in \RR^d \mid p(a) \ge 0 \text{ for every } p \in S\}
\]
where $S$ is a finite subset of $\pd:=\RR[X_1,\ldots,X_d]$. 
The main question is to compute the set
\[
\sat_>(S)=\{q \in \pd \mid q(a) > 0 \text{ for all } a \in K_S\}
\]
of all positive polynomials on $K_S$ and the set
\[
\sat_\ge(S)=\{q \in \pd \mid q(a) \ge 0 \text{ for all } a \in K_S\}
\]
of all nonnegative polynomials on $K_S$.
The answer is given by Stengle's Positivstellensatz, see \cite[Prop. 2.2.1]{mm}. 
Better answers are known under various compactness assumptions, see e.g. Jacobi's representation theorem
\cite[Theorem 5.4.4]{mm}.
We are interested in noncommutative generalizations of this theory.

Our first result, Theorem \ref{ppol}, extends the classical theory to matrix polynomials.
Fix $n \in \NN$ and write $M_n(\pd)$ (resp. $S_n(\pd)$) 
for the set of all  (resp. all symmetric) $n \times n$ matrices with entries from the set $\pd$. 
Write $\Sigma_n$ for the set of all $n \times n$ real positive semidefinite matrices. 
Let us define the positivity set of a finite subset $S \subseteq S_n(\pd)$ by
\[
K_S'=\{(a,B)  \mid a \in \RR^d, B \in \Sigma_n \setminus \{0\} \text{ and } \tr( p(a) B) \ge 0 \text{ for all } p \in S\}.
\]
Theorem \ref{ppol} computes the set
\[
\sat_>'(S)=\{q \in S_n(\pd) \mid \tr( q(a) B) > 0 \text{ for all } (a,B) \in K_S' \}
\]
under the additional assumption that $S$ contains an element of the form $(K^2-\sum_{i=1}^d X_i^2)I_n$
where $K$ is a nonzero real number and $I_n$ is the identitity matrix of size $n$. 

Note that this result is a variant of a theorem of Scherer and Hol  
(see \cite[Theorem 13]{ks} which extends \cite[Theorem 2]{hs}; a more general result follows from \cite[Theorem 3 and Lemma 5]{av}), 
where the positivity set of $S$ is defined by
\[
K_S^\mathrm{hs}=\{a \in \RR^d \mid p(a)  \text{ is positive semidefinite for all } p \in S\}
\]
and the question is to compute the set
\[
\sat_>^\mathrm{hs}(S)=\{q \in S_n(\pd) \mid q(a) \text{ is positive definite for all } a \in K_S^\mathrm{hs}\}.
\]
For $n=1$ both results reduce to Jacobi's representation theorem, see \cite[Theorem 5.4.4]{mm}. 

Our second result, Theorem \ref{tfree}, extends the classical theory to 
free $\ast$-polynomials, i.e. elements of the free algebra 
\[
F_d=\RR \langle X_1,\ldots,X_d,Y_1,\ldots,Y_d \rangle
\] 
with the involution defined by 
\[
Y_1^\ast=X_1,\ldots,Y_d^\ast=X_d.
\]
It also covers matrix versions of such polynomials, i.e.
elements of the algebra $M_n(\fd)$ with involution $[f_{ij}]^\ast=[f_{ji}^\ast]$.
Let $S$ be a finite subset of $S_n(\fd)=\{f \in M_n(\fd) \mid f^\ast=f\}$.
Its positivity set is
\[\begin{array}{ccl}
K''_S & = & \{(A_1,\ldots,A_d,B) \mid \exists m \in \NN \colon
A_1,\ldots,A_d \in M_m, B \in \Sigma_{mn} \setminus \{0\} \\
& & \quad \text{ and } \tr(f(A_1,\ldots,A_d,A_1^T,\ldots,A_d^T)B) \ge 0 \text{ for all } f \in S \}
\end{array}\]
where $M_m$ is the set of all real $m \times m$ matrices and the evaluations are performed entrywise.
Theorem \ref{tfree} computes the set 
\[\begin{array}{ccr}
\sat''_\ge(S) & = & \{f \in S_n(\fd) \mid \tr(f(A_1,\ldots,A_d,A_1^T,\ldots,A_d^T)B) \ge 0 \\
& & \text{ for all } (A_1,\ldots,A_d,B) \in K''_S\}.
\end{array}\]
For $S=\emptyset$, it extends the main theorem of \cite{mp} (which extends \cite{h}).

The closest result in literature is probably a theorem of Helton and McCullough,
see \cite[Theorem 1.2]{hm}. Here,
the positivity set of $S$ is
\[\begin{array}{ccr}
K^\mathrm{hm}_S & = & \{(A_1,\ldots,A_d) \in B(H)^d \mid f(A_1,\ldots,A_d,A_1^\ast,\ldots,A_d^\ast) \\
& & \text{ is positive definite for all } f \in S \},
\end{array}\]
where $B(H)$ is the set of all bounded operators on a real separable Hilbert space $H$
and evaluations are performed entrywise. The set
\[\begin{array}{ccl}
\sat^\mathrm{hm}_>(S) & = & \{f \in S_n(\fd) \mid f(A_1,\ldots,A_d,A_1^\ast,\ldots,A_d^\ast) \\
& & \text{ is positive definite for all } (A_1,\ldots,A_d) \in K^\mathrm{hm}_S\}
\end{array}\]
is computed under the additional assumption that $(K^2-\sum_{i=1}^d X_i^\ast X_i)I_n$ belongs to
$S$ for some nonzero real number $K$.

While the results of Scherer \& Hol and Helton \& McCullough fit into Schm\" udgen's approach
to noncommutative real algebraic geometry \cite{sch2} 
our results do not. In the last section we will formulate
an alternative approach, motivated by a paper of Helton, McCullough and Putinar \cite{hmp}, and explain its relation to our results.
The difference  is in the definition of a noncommutative real point.

\section{Matrix Polynomials}
\label{secpoly}

Let $n$ and $d$ be fixed natural numbers. Write $M_n(\pd)$ for the algebra of all
$n \times n$ matrices with entries from $\pd=\RR[X_1,\ldots,X_d]$,  
$S_n(\pd)$ for the set of symmetric matrices from $M_n(\pd)$ and $\Sigma_n(\pd)$ for the set
of all finite sums of elements of the form $P(\X)^T P(\X)$ where $P(\X) \in M_n(\pd)$.
Write also $M_n$ (resp. $S_n$, $\Sigma_n$) for the set of all (resp. all symmetric, all positive
semidefinite) $n \times n$ matrices with real entries. Let $I_n$ be the identity $n \times n$ matrix.

The aim of this section is to prove the following theorem:

\begin{thm}
\label{ppol}
Pick 
\[
P_0(\X), P_1(\X),\ldots,P_k(\X),Q(\X) \in S_n(\pd)
\]
where $P_0(\X)=(K^2-\sum_{i=1}^d X_i^2) I_n$ for some nonzero real $K$.
Suppose that for every  $a \in \RR^d$ and for every nonzero $B \in \Sigma_n$ such that
\[
\tr(P_0(a) B)\ge 0,\tr(P_1(a) B)\ge 0,\ldots,\tr(P_k(a) B)\ge 0
\] 
we have that $\tr(Q(a) B) > 0$.
Then, writing $\T$ for all sums of squares of elements from $\pd$, 
there exists $\eps \in \RR^{>0}$ such that
\[
Q(\X)-\eps I_n \in \Sigma_n(\pd)
+\T \cdot P_0(\X)+\T \cdot P_1(\X)+\ldots+\T \cdot P_k(\X). 
\]
\end{thm}

For $n=1$ we get exactly Jacobi's Representation Theorem \cite[Theorem 5.4.4]{mm}. 

\begin{proof}
We say that a subset $N$ of $S_n(\pd)$ is a weakly quadratic module (abbr. wqm) if
\[
N+N \subseteq N, \quad \T \cdot N \subseteq N \quad \mbox{ and } \quad \Sigma_n(\pd) \subseteq N.
\]
The set $\Sigma_n(\pd)$ is the smallest wqm. Clearly, the set 
\[
N_S:=\Sigma_n(\pd)+\T \cdot P_0(\X)+\T \cdot P_1(\X)+\ldots+\T \cdot P_k(\X)
\]
is the smallest wqm containing the set $S=\{P_0(\X),P_1(\X),\ldots,P_k(\X)\}$.

A wqm $N$ is said to be \textit{archimedean} if for every $A \in S_n(\pd)$ there exists 
a real $r>0$ such that $r I_n \pm A \in N$. In the sequel, we will abbreviate $I_n$ to $I$.

\subsubsection*{Step 1}
The wqm $N_S$ is archimedean.

Write $B(N_S)=\{A \in S_n(\pd) \mid  \exists r \in \RR^+ \colon r I \pm A \in N_S\}.$
Note that the set $N_S \cap \pd I$ is a quadratic module on 
$\pd I=\{pI \mid p \in \pd\}$ which contains $(K^2-\sum_{i=1}^d X_i^2)I$ for some nonzero real $K$.
By \cite[Corollary 5.2.4]{mm}, $N_S \cap \pd I$ is archimedean in $\pd I$.
It follows that $B(N_S)$ contains $\pd I$. On the other hand, it is clear that
$S_n \subseteq B(N_S)$. Suppose now that $p \in \pd$ and $A \in S_n$. 
Pick $k,l \in \RR^+$ such that $(k-p) I \in N_S$
and $l I\pm A \in N_S$. It follows $kl I \pm pA=l(k-p)I+p(l I \pm A) \in N_S$,
so that $pA \in B(N_S)$. Clearly, $B(N_S)$ is closed for addition. It follows that $S_n(\pd) 
\subseteq B(N_S)$,
by decomposing each element of $S_n(\pd)$ as a sum of products of elements from $\T$
and $S_n$. This proves the claim.

\subsubsection*{Step 2} If $M \subseteq S_n(\pd)$ is a wqm maximal subject to $-I \not\in M$
and if $M$ is archimedean, then $M \cup -M=S_n(\pd)$.

This is almost the same as \cite[Theorem 5.2.5, part (1)]{mm}. If there exists
$A \in S_n(\pd)$ such that $A \not\in M \cup -M$, then $M+\T A$ and $M-\T A$ are wqm
which strictly contain $M$. By the maximality of $M$, $-I \in M+\T A$ and $-I \in M-\T A$,
so $-I=B_1+At_1$ and $-I=B_2 -A t_2$ for some $t_1,t_2 \in \T$ and $B_1,B_2 \in M$.
Multiplying the first equality by $t_2$ and the second by $t_1$ and adding them,
we get $-(t_2+t_1)I=B_1 t_2+B_2 t_1 \in M$. It follows that $-t_1 I \in M$.
Now pick $l \in \RR^+$ such that $lI + A \in M$. It follows that
$-I=B_1+t_1(lI+ A)+l (-t_1 I) \in M+\T M+\T M \subseteq M$. This is a contradiction
with $-I \not\in M$. 

\medskip

For every wqm $N$ write $\mathcal{K}_N$ for the set of all mappings $\alpha \colon S_n(\pd) \to \RR$ such that
\begin{enumerate}
\item $\alpha(N) \ge 0$,
\item $\alpha(I)=1$,
\item $\alpha(A+B)=\alpha(A)+\alpha(B)$ for every $A,B \in S_n(\pd)$ and
\item $\alpha(tA)=\alpha(tI)\alpha(A)$ for every $t \in \pd$, $A \in S_n(\pd)$.
\end{enumerate}
For every $\alpha \in \mathcal{K}_N$, the set $\alpha^{-1}(\RR^+)$ is clearly an archimedean wqm
containing $N$.

\subsubsection*{Step 3} If $N \subseteq S_n(\pd)$ is an archimedean wqm then
$-I \not\in N$ iff $\mathcal{K}_N \ne \emptyset$.

If there is an $\alpha \in  \mathcal{K}_N$, then 
$\alpha(N) \ge 0$ and $\alpha(-I)=-1$, hence $-I \not\in N$. Conversely, if $-I \not\in N$, 
then there exists by Zorn's Lemma a wqm $M$ containing $N$ and maximal subject to $-I \not\in M$.
Since $N$ is archimedean, $M$ is also archimedean. By Step 2, $M \cup -M=S_n(\pd)$,
so we can define a mapping $\alpha \colon S_n(\pd) \to \RR$ by
\[
\alpha(A):=\sup \{r \in \RR \mid A-rI \in M\}=\inf\{r \in \RR \mid rI-A \in M\}.
\]
Clearly, $\alpha(N) \ge 0$ and $\alpha(I)=1$. The same argument as in \cite[Theorem 5.2.5, part (2)]{mm} 
shows that (3) is true for every $A,B \in S_n(\pd)$ and that
(4) is true for every $t \in \T$,$A \in M$. Since $\alpha(-A)=-\alpha(A)$ for every $A \in M$
(by (3)) and since $\T -\T=\pd$ and $M-M=S_n(\pd)$, it follows that 
(4) is true for every $t \in \pd$ and $A \in S_n(\pd)$. Thus $\alpha \in  \mathcal{K}_M$.

\subsubsection*{Step 4} A mapping $\alpha \colon S_n(\pd) \to \RR$ belongs to $\mathcal{K}_{\Sigma_n(\pd)}$
iff there exists a matrix $B \in \Sigma_n$ with $\tr(B)=1$ and a point $a \in \RR^n$ such that 
$\alpha(P(\X))=\tr(P(a)B)$ for every $P(\X) \in S_n(\pd)$.

Clearly, the homomorphism $\alpha|_{\pd I}$
is the evaluation at the point $a=(\alpha(X_1 I),\ldots,\alpha(X_d I)) \in \RR^n$
and $\alpha|_{S_n}$ is a positive functional, hence of the form $A \mapsto \tr(AB)$ for
some $B \in \Sigma_n$. Since $\alpha(I)=1$, we have that $\tr(B)=1$. 
By additivity, $\alpha(P(\X))=\tr(P(a)B)$ for every $P(\X) \in S_n(\pd)$. 
The converse is clear.

\subsubsection*{Step 5}
Suppose that $N$ is an archimedean wqm and $A=A(\X) \in S_n(\pd)$ is such that
$\alpha(A)>0$ for every $\alpha \in \mathcal{K}_N$. Then $A \in \eps I+N$
for some $\eps>0$.

The proof is the same as in \cite[Theorem 5.4.4]{mm}. Write $N_1=N-\T A$.
Since $N \subseteq N_1$, $N_1$ is archimedean. The assumption 
$\alpha(A)>0$ for every $\alpha\in \mathcal{K}_N$ implies that 
$\mathcal{K}_{N_1}=\emptyset$. By Step 3, $-I \in N_1$.
Pick $S = S(\X)\in N$ and $t \in T$ such that $-I=S-tA$, so $tA-I=S \in N$.
Since $N$ is archimedean, there exists $k \in \RR^+$ such that
$(2k-1)I-t^2A \in N$ and $(2k-t)I \in N$. Consider the identity:
$
k^2 A+(k^2r-1)I=(k-t)^2(A+rI)+2k(tA-I)+rt(2k-t)I+((2k-1)I-t^2 A).
$
This shows that $A+rI\in N$ implies that $A+(r - \frac{1}{k^2})I \in N$.
Repeating this step several times, we eventually find an $\eps >0$ such that
$A - \eps I \in N$.
\end{proof}

We can extend Theorem \ref{ppol} to the non-compact case.

\begin{thm}
Suppose that $P_1(X),\ldots,P_k(X), Q(X)\in S_n(\pd)$  are homogeneous polynomials of even degree
such that for every $a \in \RR^d \setminus \{0\}$ and for every $B \in \Sigma_n\setminus \{0\}$ which satisfy
\[
\tr(P_1(a) B)\ge 0,\ldots,\tr(P_k(a) B)\ge 0
\] 
we have that $\tr(Q(a) B) > 0$.
Then there exists $m \in \NN$ such that 
\[
(\sum_{i=1}^d X_i^2)^m Q(X) \in \Sigma_n(\pd)+\sum_{j=1}^k T P_j(X).
\]
\end{thm}

A variant of this result (corresponding to Schm\" udgen's approach) is proved in \cite{av}.
For $n=1$ both results reduce to a theorem of Putinar and Vasilescu, see  \cite{pv}.

\begin{proof}
Write $P_0(X)=(1-\sum_{i=1}^d X_i^2) \cdot 1$. From the assumption, it follows that 
$\tr(Q(a) B) > 0$ for every point $a \in \RR^n$ and for every nonzero $B \in \Sigma_n$
such that 
\[
\tr(P_0(a) B) = 0,\tr(P_1(a) B)\ge 0,\ldots,\tr(P_k(a) B)\ge 0.
\] 
By Theorem $\ref{ppol}$ there exist $q_0(X) \in \pd$, $q_1(X),\ldots,q_k(X) \in T$
and $S(X) \in \Sigma_n(\pd)$ such that
\[
Q(X) = S(X)+q_0(X)P_0(X)+q_1(X) P_1(X)+\ldots+q_k(X) P_k(X).
\]
From now on we compute in the localization of $\pd$ by $\Vert X \Vert=\sqrt{\sum X_i^2}$. 
Every elements of this localization can be written uniquely as $$\frac{g(X)+h(X)\Vert X \Vert}{\Vert X \Vert^l} $$
where $g(X),h(X) \in \pd$  and $l \in \NN$.
Since $P_0(\frac{X}{\Vert X \Vert})=0$, we get
\[
Q(\frac{X}{\Vert X \Vert})=S(\frac{X}{\Vert X \Vert})+\sum_{j=1}^k q_j(\frac{X}{\Vert X \Vert})P_j(\frac{X}{\Vert X \Vert}).
\]
By clearing denominators and comparing components at $1$ and $\Vert X \Vert$, we get 
(because $Q(X)$ and $P_j(X)$ are homogeneous of even degree) that 
\[
\Vert X \Vert^{2m} Q(X) \in \Sigma_n(\pd)+\sum_{j=1}^k T P_j(X)
\]
for some $m \in \NN$.
\end{proof}

\section{Matrix free $\ast$-polynomials}
\label{secfree}

Let $n$ and $d$ be fixed natural numbers. Write $\fd$ for the free $\ast$-algebra
in $d$ variables as defined in the introduction and $M_n(\fd)$ for the $\ast$-algebra
of $n \times n$ matrices with entries from $\fd$ with involution defined by $[f_{ij}]^\ast=[f_{ji}^\ast]$.
Write $S_n(\fd)$ for the real vector space of all elements $f \in M_n(\fd)$ such that $f^\ast=f$
and $\Sigma_n(\fd)$ for the set of all finite sums of elements of the form $f^\ast f$, $f \in M_n(\fd)$.
Let $M_n$, $S_n$ and $\Sigma_n$ be as in the previous section.

For every $d$-tuple $\mathbf{C}=(C_1,\ldots,C_d) \in (M_m)^d$, we define 
the mapping $\ev_\mathbf{C} \colon \fd \to M_{m}$
by
\[
\ev_\mathbf{C}(f(X_1,\ldots,X_d,X_1^\ast,\ldots,X_d^\ast))
=f(C_1,\ldots,C_d,C_1^T,\ldots,C_d^T)
\]
and the mapping $(\ev_\mathbf{C})_n \colon M_n(\fd) \to M_{mn}$
\[
(\ev_\mathbf{C})_n([f_{ij}])=[\ev_\mathbf{C}(f_{ij})].
\]

The aim of this section is to prove the following theorem.

\begin{thm}
\label{tfree}
Pick $n$ and $d$. For every elements $p_1,\ldots,p_k,q \in S_n(\fd)$,
the following are equivalent:
\begin{enumerate}
\item $\tr((\ev_\mathbf{C})_n(q)B) \ge 0$ for every $m \in \NN$, $\mathbf{C} \in (M_m)^d$ and $B \in \Sigma_{mn}$
such that $\tr((\ev_\mathbf{C})_n(p_i)B) \ge 0$ for every $i=1,\ldots,k$.
\item $q \in \overline{\Sigma_n(\fd)+\RR^+ p_1+\ldots+\RR^+ p_k}$, where the closure refers 
to the finest locally convex vector space topology of $S_n(\fd)$.
\end{enumerate}
\end{thm}

The proof will depend on some rather elementary observations about $\ast$-representations that we are going to recall now.
Recall that a \textit{$\ast$-representation} of $\A$ is an ordered pair $(\pi,D_\pi)$
where $D_\pi$ is a real inner product space and $\pi$ is a unital real algebra 
homomorphism from $\A$ into the algebra of all linear operators on $D_\pi$
such that $\langle \pi(a) v_1,v_2 \rangle= \langle v_1, \pi(a^\ast) v_2 \rangle$
for every $a \in \A$. 

\begin{ex}
For every $m \in \NN$ and $\mathbf{C} \in (M_m)^d$, the mapping $\ev_\mathbf{C} \colon \fd \to M_{m}$
defines a $\ast$-representation $(\ev_\mathbf{C}, \RR^m)$ of the $\ast$-algebra $F_d$.
Conversely, every $\ast$-representation $(\pi,D_\pi)$ of $F_d$ for which $m=\dim D_\pi < \infty$
comes from $\ev_\mathbf{C}$ with $\mathbf{C}=(\pi(X_1),\ldots,\pi(X_d))\in (M_m)^d$.
\end{ex}

Let $\A$ be a $\ast$-algebra and $(\pi,D_\pi)$ a $\ast$-representation of $\A$.
We can equip $D_\pi$ with the structure of a left $\A$-module 
by setting $av:=\pi(a)v$ for every $a \in \A$ and $v \in D_\pi$. Equivalently,
we can start with a left $\A$-module $D$ equipped with a real valued
inner product satisfying $\langle a  v_1,v_2 \rangle
=\langle v_1,a^\ast  v_2 \rangle$ for every $a \in \A$ and $v_1,v_2 \in D$
and denote its action by $\pi$. We say that the $\ast$-representations
$(\pi,D_\pi)$ and $(\psi,D_\psi)$ are \textit{unitarily equivalent} 
if there exist mutually inverse isometries $S \colon D_\pi \to D_\psi$
and $T \colon D_\psi \to D_\pi$ such that $S \pi(a)=\psi(a) S$ and
$T \psi(a)=\pi(a)T$ for every $a \in \A$. Equivalently, $S$ and $T$ are mutually inverse 
left $\A$-module homomorphisms between $D_\pi$ and $D_\psi$ that preserve inner products.

\begin{lemma}\label{l1}(GNS construction) 
Let $\A$ be a $\ast$-algebra, $\Ah=\{a \in \A \mid a=a^\ast\}$ and $\Sigma_\A$
the set of all finite sums of elements $a^\ast a$, $a \in \A$.
Let $f$ be a real linear functional on $\Ah$ such that $f(\Sigma_\A) \ge 0$. Then there exists a
$\ast$-representation $(\pi_f,D_f)$ of $\A$ and a vector $v_f \in D_f$
such that $f(a)=\langle \pi_f(a)v_f,v_f \rangle$ for every $a \in \Ah$.
\end{lemma}

\begin{proof}
The set $I_f=\{a \in \A \mid f(a^\ast a)=0\}$ is a left ideal in $\A$
and $D_f:=\A/I_f$ is a left $\A$-module with inner product 
$\langle x+I_f,y+I_f \rangle:=\frac12 f(y^\ast x+x^\ast y)$. Let $\pi_f$ be the action of $D_f$
and $v_f:=1+I_f \in D_f$. See Theorem 8.6.2 in \cite{sch} for details.
\end{proof}

\begin{lemma}\label{l2}(Morita equivalence)
Let $n$ be a natural number, $\A$ a $\ast$-algebra and $M_n(\A)$ the $\ast$-algebra
of $n \times n$ matrices with entries in $\A$ with involution $[a_{ij}]^\ast:=[a_{ji}^\ast]$. 
Every $\ast$-representation $(\pi,D_\pi)$ of $\A$ induces a $\ast$-representation
$(\pi_n ,D_{\pi_n})$ of $M_n(\A)$ where $D_{\pi_n}:=(D_\pi)^n$ and
$\pi_n([a_{ij}]):=[\pi(a_{ij})]$. Moreover, every $\ast$-representation
of $M_n(\A)$ is unitarily equivalent to a $\ast$-representation induced from $\A$ as above.
\end{lemma}

\begin{proof}
The first part is clear. 
Write $E_{ij}$ for the element of $M_n(\A)$ which has $1$ at $(i,j)$-th place and zeros elsewhere.
Pick a $\ast$-representation $(\psi,D_\psi)$ of $M_n(\A)$ and note that $D:=E_{11} D_\psi$
is a left $\A$-submodule of $D_\psi$. We equip $D$ with the inner product inherited from $D_\psi$. 
Write $\pi$ be the action of $D$ and note that $D_{\pi_n}=D^n$.
The mappings $v \mapsto (E_{11}v,\ldots,E_{1n}v)$ from $D_\psi$ to $D_{\pi_n}$ and 
$(v_1,\ldots,v_n) \mapsto E_{11}v_1+\ldots+E_{n1}v_n$ from $D_{\pi_n}$ to $D_\psi$
are mutually inverse homomorphisms of left $M_n(\A)$-modules which preserve inner products.
\end{proof}

The following is similar to Proposition 4 in \cite{sch2}. For every element of $\fd$
we can define its \textit{degree} as the total degree in $X_i$ and $X_j^\ast$. 
Write $\fdk$ for the set of all elements of $\fd$ of degree $\le K$.

\begin{lemma}
\label{l3}
For every $\ast$-representation $(\pi,D_\pi)$ of $\fd$, every
finite subset $\{v_1,\ldots,v_n\}$ of $D_\pi$ and every $K \in \NN$ 
there exists a $\ast$-representation $(\rho,V)$ of $\fd$ where $V
= \pi(\fdk)v_1+\ldots+\pi(\fdk)v_k$ is a finite-dimensional subspace 
of $D_\pi$ containing $\{v_1,\ldots,v_n\}$
and $\rho(b) v_i=\pi(b) v_i$ for every $b \in \fdk$ and $i=1,\ldots, n$.
\end{lemma}

\begin{proof}
Let $H$ be the Hilbert space completion of $D_\pi$ and let $P$ be the orthogonal projection
of $H$ to $V$. For every $x \in \{X_1,\ldots,X_d,X_1^\ast,\ldots,X_d^\ast\}$ write $\rho(x) = P \pi(x)|_V$.
For every $u,u' \in V$, we have that $\langle \rho(x) u,u' \rangle_V=\langle \rho(x) u,u' \rangle_H=\langle P\pi(x) u,u' \rangle_H
=\langle \pi(x) u,P^\ast u' \rangle_H=\langle \pi(x) u,u' \rangle_H=
\linebreak 
\langle u, \pi(x)^\ast u' \rangle_H= 
\langle u, \pi(x^\ast) u' \rangle_H=\langle P^\ast u, \pi(x^\ast) u' \rangle_H 
=\langle u, P \pi(x^\ast) u' \rangle_H =\langle u, \rho(x^\ast) u' \rangle_H =\langle u, \rho(x^\ast) u' \rangle_V$,
hence $\rho(x^\ast)=\rho(x)^\ast$. Since $\fd$ is a free $\ast$-algebra,
we can extend $\rho$ to a $\ast$-representation of $\fd$.

Fix $i$ between $1$ and $n$.
To prove that $\rho(b) v_i=\pi(b) v_i$ for every $b \in \fdk$ we may assume that
$b$ is a monomial in $X_i$ and $X_j^\ast$ and proceed by induction on $\deg b$.
If $b=1$, this is clear. Suppose that $b=xc$ where $b,c \in \fdk$ and $x \in \{X_1,\ldots,X_d,X_1^\ast,\ldots,X_d^\ast\}$.
By the inductive hypothesis, we have that $\rho(c)v_i=\pi(c)v_i$. It follows that
$\rho(b)v_i=\rho(x)\rho(c)v_i=\rho(x)\pi(c)v_i$. Since  $\pi(c)v_i \in V$, 
we have that $\rho(x)\pi(c)v_i=P \pi(x) \pi(c)v_i$ by the definition of $\rho(x)$. Since
$b \in \fdk$, we have that $\pi(b)v_i \in V$, hence $P \pi(x) \pi(c)v_i=P \pi(b)v_i=\pi(b)v_i$.
Therefore, $\rho(b) v_i=\pi(b) v_i$.
\end{proof}

We are now able to give the proof of Theorem \ref{tfree}.

\begin{proof}
Recall that the finest locally convex topology of a real vector space is the topology whose fundamental
system of neighbourhoods of zero consists of all convex absorbing sets. In particular, every
real linear functional on $S_n(\fd)$ is continuous in this topology. By the Separation theorem, 
an element $q$ of $S_n(F_d)$ belongs to the closure of a convex cone $C \subseteq S_n(F_d)$ if and only 
if $f(q) \ge 0$ for every real linear functional $f$ on $S_n(F_d)$ such that $f(C) \ge 0$.
For $C=\Sigma_n(\fd)+\RR^+ p_1+\ldots+\RR^+ p_k$, we get that (2) is equivalent to
\begin{enumerate}
\item[(A)] $f(q) \ge 0$ for every real linear functional $f$ on $S_n(\fd)$ such that
$f(\Sigma_n(\fd)) \ge 0$ and $f(p_1) \ge 0, \ldots,f(p_k) \ge 0$.
\end{enumerate}
For every $\ast$-representation $\pi$ of $\A$ and for every $v \in D_\pi$,
the real linear functional $f_{(\pi,v)}(a):=\langle \pi(a)v,v \rangle$, $a \in S_n(\fd)$,
satisfies $f_{(\pi,v)}(\Sigma_n(\fd))\ge 0$. Conversely, for every real linear functional $f$
on $S_n(\fd)$ such that $f(\Sigma_n(\fd)) \ge 0$ there exist by Lemma \ref{l1} 
a $\ast$-representation $(\pi_f,D_f)$ of $M_n(\fd)$ and a vector $v_f \in D_f$ such that
$f(a)=\langle \pi_f(a)v_f,v_f \rangle$ for every $a \in S_n(\fd)$. It follows that (A) is equivalent to:
\begin{enumerate}
\item[(B)] $\langle \pi(q)v,v \rangle \ge 0$ for every $\ast$-representation $(\pi,D_\pi)$ of $M_n(\fd)$ and every $v \in D_\pi$ 
such that $\langle \pi(p_1)v,v \rangle \ge 0, \ldots,\langle \pi(p_k)v,v \rangle \ge 0$.
\end{enumerate}
By Lemma \ref{l2} every $\ast$-representation of $M_n(\fd)$ is
unitarily equivalent to a $\ast$-representation of the form $\psi_n$. Therefore, (B) is equivalent to
\begin{enumerate}
\item[(C)] $\langle \psi_n(q)u,u \rangle \ge 0$ for every $\ast$-representation $(\psi,D_\psi)$ of $\fd$ and every 
$u \in (D_\psi)^n$ such that $\langle \psi_n(p_i)u,u \rangle \ge 0$ for all $i=1,\ldots,k$.
\end{enumerate}
For every $\ast$-representation $(\psi,D_\psi)$ of $\fd$ and every $u=(u_1,\ldots,u_n) \in (D_\psi)^n$, 
there exists by Lemma \ref{l3} a finite-dimensional $\ast$-representation $(\rho,D_\rho)$ such that $\phi(a)u_i=\rho(a)u_i$
for every $a \in \fdk$ and every $i=1,\ldots,n$. It follows that $\langle \psi_n(a)u,u \rangle
= \langle \rho_n(a)u,u \rangle$ for every $a \in \fd$ and every $u \in (D_\psi)^n$,
hence (C) is equivalent to
\begin{enumerate}
\item[(D)] $\langle \rho_n(q)u,u \rangle \ge 0$ for every finite-dimensional $\ast$-representation $(\rho,D_\rho)$ of $\fd$ and every 
$u \in (D_\psi)^n$ such that $\langle \rho_n(p_i)u,u \rangle \ge 0$ for all $i=1,\ldots,k$.
\end{enumerate}
If $(\rho,D_\rho)$ is a finite-dimensional $\ast$-representation of $\fd$ and $C_1, \ldots, C_d$ are matrices that
belong to the operators $\rho(X_1),\ldots,\rho(X_d)$ in some orthonormal basis of $D_\rho$, then $(\rho,D_\rho)$
is unitarily equivalent to the $\ast$-representation $(\ev_\mathbf{C},\RR^m)$, where $\mathbf{C}=(C_1,\ldots,C_d)$,
and $(\rho_n,(D_\rho)^n)$ is unitarily equivalent to $((\ev_\mathbf{C})_n,\RR^{mn})$. Hence, (D) is equivalent to
\begin{enumerate}
\item[(E)] $\langle (\ev_\mathbf{C})_n(q)v,v \rangle \ge 0$ for every $m \in \NN$, every  $\mathbf{C} \in (M_m)^d$
and every $v \in \RR^{mn}$ such that  $\langle (\ev_\mathbf{C})_n(p_i)v,v \rangle \ge 0$ for all $i$.
\end{enumerate}
For every $d$-tuple $\mathbf{C}=(C_1,\ldots,C_d) \in (M_m)^d$, 
write $mn \oplus \mathbf{C}$ for the $d$-tuple $(mn \oplus C_1,\ldots,mn \oplus C_d)
\in (M_{m^2 n})^d$. Substituting $(\mathbf{C},v) \to (mn \oplus \mathbf{C},(u_1,\ldots,u_{mn}))$ in one direction and
$(mn \oplus \mathbf{C},(u_1,\ldots,u_{mn})) \to (mn \oplus \mathbf{C},v \otimes e_1)$ 
in the other direction, where $u_1,\ldots,u_{mn} \in \RR^{mn}$, $e_1=(1,0,\ldots,0) \in \RR^{mn}$
and $\otimes$ is the Kronecker product, we see that (E) is equivalent to
\begin{enumerate}
\item[(F)] $\langle (\ev_{mn \oplus \mathbf{C}})_n(q)(u_1,\ldots,u_{mn}),(u_1,\ldots,u_{mn}) \rangle \ge 0$ for every 
$m \in \NN$, every  $\mathbf{C} \in (M_m)^d$ and every $(u_1,\ldots,u_{mn}) \in (\RR^m)^{mn}$ such that 
$\langle \ev_{mn \oplus \mathbf{C}}(p_i)(u_1,\ldots,u_{mn}),(u_1,\ldots,u_{mn})\rangle \ge 0$ for all $i$.
\end{enumerate}
We claim that (F) is equivalent to (1).
Let $P$ be the matrix of the transpose mapping on $M_{mn}$ in the basis $E_{11},\ldots,E_{1n},\ldots,E_{n1},\ldots,E_{nn}$.
Note that $P$ is a permutation (hence orthogonal) matrix such that $B \otimes I_{mn} = P^T (I_{mn} \otimes B)P$ for every 
$B \in M_{mn}$. In particular, we have that $\ev_{mn \oplus \mathbf{C}}(f)=
\ev_\mathbf{C}(f) \otimes I_{mn}=P^T(I_{mn} \otimes \ev_\mathbf{C}(f))P=P^T(mn \oplus \ev_\mathbf{C}(f))P$ for every 
$f \in M_n(\fd)$. Writing $u=(u_1,\ldots,u_{mn})$ and $v=Pu=(v_1,\ldots,v_{mn})$, we get that
$
\langle (\ev_{mn \oplus \mathbf{C}})_n(q)u,u \rangle 
= \langle P^T(mn \oplus (\ev_{\mathbf{C}})_n(q))Pu,u \rangle
= \langle (mn \oplus (\ev_{\mathbf{C}})_n(q))v,v \rangle
= \sum_{i=1}^{mn} v_i^T (\ev_{\mathbf{C}})_n(q)v_i 
= \tr((\ev_{\mathbf{C}})_n(q)B)
$
where $B=\sum_{i=1}^{mn} v_i v_i^T$ belongs to $\Sigma_{mn}$.
Conversely, every element $B$ of $\Sigma_{mn}$ has a decomposition $B=\sum_{i=1}^{mn} v_i v_i^T$ with $v_i \in \RR^{mn}$ and 
for $u=P^T(v_1,\ldots,v_n)$ we have that $\langle (\ev_{mn \oplus \mathbf{C}})_n(q)u,u \rangle = \tr((\ev_{\mathbf{C}})_n(q)B)$.
\end{proof}

A subset $C$ of a real vector space $V$ is a \textit{convex cone} if $C+C \subseteq C$ and $\RR^+ C \subseteq C$.
Write $\overline{C}$ for the closure and $C^\circ$ for the interior of $C$ in the finest locally convex topology of $V$.
Write $C^\vee$ for the set of all real linear functionals $f$ on $V$
such that $f(C) \ge 0$. The following is well known and easy to prove. 

\begin{prop}
\label{qprop}
Let $C$ be a convex cone in real vector space $V$. An element $\q \in V$ belongs to $C^\circ$
if and only if for every element $v \in V$ there exists $n \in \NN$ such that $n\q+v \in C$.
If one of these equivalent conditions is true then for every $v \in V$ the following are equivalent:
\begin{enumerate}
\item[(1)] $v \in \overline{C}$, 
\item[(2)] $v+\eps \q \in C$ for every real $\eps >0$, and 
\item[(3)] $f(v) \ge 0$ for every $f \in C^{\vee}$. 
\end{enumerate}
Under the same assumptions, the following are equivalent:
\begin{enumerate}
\item[(1')] $v \in C^\circ$, 
\item[(2')] $v-\eps \q \in C$ for some real $\eps>0$, and 
\item[(3')] $f(v)>0$ for every $f \in C^\circ$.
\end{enumerate}
\end{prop}

The following is the archimedean version of Theorem \ref{tfree}.

\begin{thm}
\label{archfree}
Pick $n$ and $d$. For every elements $p_1,\ldots,p_k\in S_n(\fd)$
such that the cone $C:=\Sigma_n(\fd)+\RR^+ p_1+\ldots+\RR^+ p_k$
is archimedean (i.e. $1 \in C^\circ$) the following are equivalent 
for every $q \in S_n(\fd)$:
\begin{enumerate}
\item[(a)] $\tr((\ev_\mathbf{C})_n(q)B) \ge 0$ for every $m \in \NN$, 
$\mathbf{C} \in (M_m)^d$ and $B \in \Sigma_{mn}$ such that 
$\tr((\ev_\mathbf{C})_n(p_i)B) \ge 0$ for all $i=1,\ldots,k$.
\item[(b)] $q + \eps I_n \in C$ for every real $\eps >0$.
\end{enumerate}
Under the same assumptions, the following are equivalent:
\begin{enumerate}
\item[(a')] $\tr((\ev_\mathbf{C})_n(q)B) > 0$ for every $m \in \NN$,  
$\mathbf{C} \in (M_m)^d$ and  $B \in \Sigma_{mn}\setminus\{0\}$ such that 
$\tr((\ev_\mathbf{C})_n(p_i)B) \ge 0$ for all $i=1,\ldots,k$.
\item[(b')] $q - \eps I_n \in C$ for some real $\eps >0$.
\end{enumerate}
\end{thm}

\begin{proof}
By Theorem \ref{tfree}, (a) is equivalent to $q \in \overline{C}$
and by Proposition \ref{qprop} applied to $V=S_n(\fd)$ and $e=I_n$,
$q \in \overline{C}$ is equivalent to (b).

To prove the second part, consider the following claim:
\begin{enumerate}
\item[(A')] $f(q) > 0$ for every real linear functional $f$ on $S_n(\fd)$ such that
$f(\Sigma_n(\fd)) \ge 0$ and $f(p_1) \ge 0, \ldots,f(p_k) \ge 0$.
\end{enumerate}
We can prove that (A') is equivalent to (a') by following the proof
of Theorem \ref{tfree} (the equivalence of assertions (A) and (1)).
Applying Proposition \ref{qprop} as above, we see that (A') is equivalent to (b').
\end{proof}

\section{General $\ast$-algebras}

Let $\A$ be a fixed $\ast$-algebra and $\rall=\rall(\A)$ the class of all $\ast$-representations of $\A$.
In K. Schm\" udgen's approach to noncommutative real algebraic the elements of 
$\Ah=\{a \in \A \mid a=a^\ast\}$ are considered as `noncommutative real polynomials'
and the elements of some fixed subclass $\R$ of $\rall$ are considered as `noncommutative real points'.
Interesting choices for $\R$ include the class $\rfin$ of all finite-dimensional $\ast$-representations
of $\A$ and the class $\rbd$ of all bounded $\ast$-representations of $\A$.
The \textit{positivity set} of a given subset $S$ of $\Ah$ is defined by
\[
K_S^{\R,\sch}:=\{ (\pi,D_\pi) \in \R \mid \pi(s) \succeq 0 \text{ for all } s \in S \}
\]
and the problem is to compute the \textit{saturations}
\[
\sat^{\R,\sch}_>(S) := \{q \in \Ah \mid \pi(q)\succ 0 \text{ for all } (\pi,D_\pi) \in K^{\R,\sch}_S \}
\]
\[
\sat^{\R,\sch}_\ge(S) := \{q \in \Ah \mid \pi(q)\succeq 0 \text{ for all } (\pi,D_\pi) \in K^{\R,\sch}_S \}
\]
and
\[
\sat^{\sch,\R}_{\not\le}(S):=\{a \in \sym(\A) \mid \pi(a) \not\preceq 0 \text{ for all } (\pi,D_\pi) \in K_S^{\sch,\R}\}.
\]

\begin{ex}
\label{esch}
If $\A=M_n(\pd)$, $S \subseteq \Ah$ and $\R(\A)$ is the set of all  mappings
$\ev_a \colon \A \to M_n$, $\ev_a(P(X))=P(a)$, where $a \in \RR^d$, then
\[
\sat^\mathrm{hs}_>(S)=\sat^{\R(\A),\sch}_>(S).
\]
If  $\B =M_n(\fd)$ and $S \subseteq \Bh$ then
\[
\sat^\mathrm{hm}_>(S)=\sat^{\rbd(\B),\sch}_>(S).
\]
The sets $\sat^\mathrm{hs}_>(S)$ and $\sat^\mathrm{hm}_>(S)$ were defined in the introduction.
\end{ex}

Recall that a subset $M$ of $\Ah$ is called a \textit{quadratic module} in $\A$ 
if $M+M \subseteq M$, $1 \in M$ and $a^\ast M a \subseteq M$ for every $a \in \A$.
Write $\overline{M}$ for the closure and $M^\circ$ for the interior of
$M$ in the finest locally convex topology of the real vector space $\Ah$.
A quadratic module $M$ is \textit{archimedean} if for every $a \in \Ah$ there
exists a real positive $k$ such that $k \cdot 1+a \in M$. By Proposition \ref{qprop},
$M$ is archimedean iff $1 \in M^\circ$.

\begin{thm}
\label{tsch}
Let $\A$ be a $\ast$-algebra and $S \subseteq \Ah$. Write $M_S$ for the smallest quadratic module
in $\A$ which contains $S$. Then
\[
\sat^{\rall,\sch}_\ge(S)=\overline{M_S}
\]
Moreover, if $M_S$ is archimedean then
\[
\sat^{\rbd,\sch}_>(S) =\{a \in \Ah \mid \eps \cdot 1+a \in M_S \text{ for some } \eps \in \RR^{>0} \}= (M_S)^\circ,
\]
\[
 \sat^{\rbd,\sch}_\ge(S)  =\{a \in \Ah \mid \eps \cdot 1+a \in M_S \text{ for all } \eps \in \RR^{>0} \} = \overline{M_S},
\]
\[
\quad \sat^{\rbd,\sch}_{\not\le}(S)= \{f \in \sym(\A) \mid -1 \in M_{S \cup \{-f\}} \, \}.
\]
\end{thm}

The proof of the first part is almost the same as the proof of (i) $\Leftrightarrow$ (ii)
in \cite[Proposition 3]{sch2} or the proof of  (2) $\Leftrightarrow$ (B) in our Theorem \ref{tfree}.
The second part is the same as \cite[Theorem 12]{c1} and \cite[Theorem 5]{c2}, or
Propositions 14-16 in \cite{sch2}.
Note that the theorem of  Helton \& McCullough follows from Theorem \ref{tsch} and Example \ref{esch}.
For  Scherer \& Hol you also need an observation from the proof of \cite[Theorem 13]{ks}
that $\sat^{\R(\A),\sch}_>(S)=\sat^{\rbd(\A),\sch}_>(S)$ when $M_S$ is archimedean.

Our main results do not fit into Schm\" udgen's approach but they fit instead into the approach 
that was outlined (in the case of free algebras) by Helton, McCullough and Putinar in \cite{hmp}. 
In this approch `noncommutative real polynomials' are the same as above, i.e. elements of $\Ah$,
but `noncommutative real points' are different - they are triples
$(\pi,D_\pi,v)$, where $(\pi,D_\pi)$ belongs to $\R$ and $v$ belongs to $D_\pi \setminus \{0\}$.
Write $\pts$ for the set of all such triples.

Let $S$ be a subset of $\Ah$. Its positivity set is defined by
\[
K_S^{\R,\hmp} := \{ (\pi,D_\pi,v) \in \pts \mid 
\langle \pi(s) v,v \rangle \ge 0 \text{ for all } s \in S \}
\]
and the corresponding saturations are defined by
\[
\sat^{\R,\hmp}_>(S) := \{q \in \Ah \mid \langle \pi(q)v,v \rangle > 0 \text{ for all } (\pi,D_\pi,v) \in K^{\R, \hmp}_S\}
\]
and
\[
\sat^{\R,\hmp}_\ge(S) := \{q \in \Ah \mid \langle \pi(q)v,v \rangle \ge 0 \text{ for all } (\pi,D_\pi,v) \in K^{\R, \hmp}_S\}.
\]

\begin{ex}
Recall the definition of the set $\sat_{\ge}'(S)$ from the introduction.
If $\A=M_n(\pd)$, $S \subseteq \Ah$ and  $\R'(\A)=\{n \oplus \ev_a \mid a \in \RR^d\}$
where $(n \oplus \ev_a)(P) = n \oplus P(a)$ for every $P=P(X) \in \A$, then
\[
\sat_{\ge}'(S)=\sat^{\R'(\A),\hmp}_{\ge}(S).
\]
This follows from the identity
$
\langle (n \oplus \ev_a)(P)(v_1,\ldots,v_n), (v_1,\ldots,v_n) \rangle 
= \sum_{i=1}^{n} v_i^T \ev_a(P)v_i 
= \tr(\ev_a(P)B)
$
where $B=\sum_{i=1}^{n} v_i v_i^T \in \Sigma_{n}$.
(Compare with equivalence (F) $\Leftrightarrow$ (1) in the proof of Theorem \ref{tfree}.)

We would like to point out that in general
\[
\sat^{\R(\A),\hmp}_{\ge}(S) \supsetneq \sat^{\R'(\A),\hmp}_{\ge}(S) \supsetneq
\sat^{\rfin(\A),\hmp}_{\ge}(S) 
\]
Both inclusions are clear
(because larger representation class means larger positivity set and so smaller saturation);
we just have to prove that they are proper. 

(a) For the following constant elements of $S_2(\pd)$
\[
p_1=\left[ \begin{array}{cc} 1 & 0 \\ 0 & -1 \end{array} \right],
p_2=\left[ \begin{array}{cc} 0 & -1 \\ -1 & 1 \end{array} \right]
 \mbox{   and   } 
Q=\left[ \begin{array}{cc} 0 & -1 \\ -1 & -1 \end{array} \right]
\]
we have that
\[
Q \in \sat^{\rev(\A),\hmp}_{\ge}(\{p_1,p_2\})\setminus\sat^{\R(\A),\hmp}_{\ge}(\{p_1,p_2\}).
\]
\textit{Proof.} We claim that for every $v=[x,y]^T \in \RR^2$ such that
$v^\ast p_1 v \ge 0$ and $v^\ast p_2 v \ge 0$.
we also have that $v^\ast Q v \ge 0$. The claim is clearly true if $x=0$.
If $x\ne 0$, then we can assume that $x=1$. We get 
$v^\ast p_1 v \ge 0$ iff $y^2 \le 1$,
$v^\ast p_2 v \ge 0$ iff $(y-1)^2 \ge 1$ and
$v^\ast Q v \ge 0$ iff $(y+1)^2 \le 1$, which implies the claim.
On the other hand, note that $e_1^T p_1 e_1+e_2^T p_1 e_2 \ge 0$
and $e_1^T p_2 e_1+e_2^T p_2 e_2 \ge 0$ but $e_1^T Q e_1+e_2^T Q e_2 < 0$ 
where $e_1,e_2$ is the standard basis of $\RR^2$.

(b) For 
$
p=X_1 \mbox{  and  } q=X_1^3
$
from  $\A=\RR[X_1]$ we have that
\[
q \in \sat^{\R(\A),\hmp}_{\ge}(\{p\})\setminus\sat^{\rfin(\A),\hmp}_{\ge}(\{p\}).
\]
\textit{Proof.} Since $a_1 \ge 0$ implies that $a_1^3 \ge 0$ for every $a_1 \in \RR$,
it follows that $q \in \sat^{\R(\A),\hmp}_{\ge}(\{p\})$. On the other hand, for
$\pi \colon u(X) \mapsto [u(-2)] \oplus [u(1)]$ and $v=[1,2]$, we have that 
$\langle\pi(p(X))v,v \rangle =p(-2)1^2+p(1)2^2=2 \ge 0$ and
$\langle\pi(q(X))v,v \rangle =q(-2)1^2+q(1)2^2=-4 < 0$, hence
$q \not\in \sat^{\rfin(\A),\hmp}_{\ge}(\{p\})$.
\end{ex}

\begin{ex} Recall the definition of the set $\sat_{\ge}''(S)$ from the introduction. 
If  $\B =M_n(\fd)$, $S \subseteq \Bh$ and $\R(\B)$ is the set of all mappings $(\ev_\mathbf{C})_n$
where $\mathbf{C}$ is a $d$-tuple of same size matrices then
\[
\sat_{\ge}''(S)=\sat^{\R(\B),\hmp}_{\ge}(S).
\]
By the proof of Theorem \ref{tfree} we also have that
\[
\sat^{\R'(\B),\hmp}_{\ge}(S)=\sat^{\rfin,\hmp}_{\ge}(S)=\sat^{\rall(\B),\hmp}_{\ge}(S).
\]
\end{ex}

A subset $N$ of $\Ah$ will be called a \textit{quadratic cone} in $\A$
if $N+N \subseteq N$, $\RR^{\ge 0} \cdot N \subseteq N$ and $a^\ast a \in N$
for every $a \in \A$. 
The following general result can be extracted from the proofs
of Theorems \ref{tfree} and \ref{archfree}.
It is an analogue of Theorem \ref{tsch}.

\begin{thm}
\label{thmp}
Let $\A$ be a $\ast$-algebra and $S \subseteq \Ah$. Write $N_S$ for the smallest quadratic cone
in $\A$ which contains $S$. Then
\[
\sat^{\rall,\hmp}_\ge(S)=\overline{N_S}.
\]
Moreover, if $N_S$ is archimedean then
\[
\sat^{\rbd,\hmp}_\ge(S)  =\{a \in \Ah \mid \eps \cdot 1+a \in N_S \text{ for all } \eps \in \RR^{>0} \} =\overline{N_S}
\]
and
\[
\sat^{\rbd,\hmp}_>(S) =\{a \in \Ah \mid \eps \cdot 1+a \in N_S \text{ for some } \eps \in \RR^{>0} \}= (N_S)^\circ.
\]
\end{thm}

For $\A=M_n(\fd)$, Theorem \ref{thmp} extends Theorem \ref{tfree}. However, if $\A=M_n(\pd)$, it extends
neither Theorem \ref{ppol} nor the solution of (the usual or the matrix) Hilbert's 17th problem. The matrix
version of the Hilbert's 17th problem was solved independently by \cite{gr} and \cite{ps}. For a constructive
proof see Proposition 10 in \cite{sch} which precedes \cite{hn}.

For $S = \emptyset$, Theorem \ref{thmp} is the same as Theorem \ref{tsch}. The case $N_\emptyset$ archimedean
(i.e. $\A$ algebraically bounded) is known as Vidav-Handelman theory, see \cite[Section 1]{han} and \cite{vidav}.

\end{document}